\documentclass{amsart}

\usepackage{amssymb}

\newtheorem{theorem}{Theorem}[section]

\newtheorem{claim}[theorem]{Claim}

\numberwithin{equation}{section}

\begin{document}
\setcounter{page}{1}

\title[Non-intersecting random walks]
{Non-intersecting, simple, symmetric random walks and the extended
Hahn kernel}
\author[K.~Johansson]{Kurt Johansson}

\address{
Department of Mathematics,
Royal Institute of Technology,
S-100 44 Stockholm, Sweden}

\email{kurtj@math.kth.se}

\thanks{Supported by the Swedish Science Research Council and the
G\"oran Gustafsson Foundation (KVA)}



\dedicatory{To Pierre van Moerbeke on his 60:th birthday}

\begin{abstract} Consider $a$ particles performing simple, symmetric,
non-intersecting random walks, starting at points $2(j-1)$, $1\le j\le a$
at time 0 and ending at $2(j-1)+c-b$ at time $b+c$. This can also be 
interpreted as a random rhombus tiling of an $abc$-hexagon, or as a
random boxed planar partition confined to a rectangular box with side
lengths $a$, $b$ and $c$. The positions of the particles at all times gives a
determinantal point process with a correlation kernel given in terms
of the associated Hahn polynomials. In a suitable scaling limit
we obtain non-intersecting Brownian motions which can be related to 
Dysons's Hermitian Brownian motion via a suitable transformation.
\end{abstract}

\maketitle

\section{Introduction}

We will consider $a$ simple, symmetric random walks started at $2(j-1)$,
$1\le j\le a$, conditioned not to intersect in the time interval
$[0,b+c]$, and end at $c-b+2(j-1)$ at time $b+c$. Here $a,b,c$, $c\ge b$, 
fixed positive integers. This model has several interpretations. One is 
as a uniform random rhombus tiling of an $abc$-hexagon, i.e. a hexagon with
side lengths $a,b,c,a,b,c$, see \cite{CLP}. This translates directly to  
a dimer or
perfect matching representation, see e.g. \cite{Ke}, so it is a 
kind of two-dimensional statistical mechanics model.
Another interpretation is as a boxed planar partition in a rectangular
box with side lengths $a$, $b$ and $c$, \cite{Sta}. The number of 
possible configurations, the partition function of the model, was 
computed by MacMahon, and is given by
\begin{equation}\label{1.1}
Z(a,b,c)=\prod_{i=1}^a\prod_{j=1}^b\prod_{k=1}^c\frac{i+j+k-1}
{i+j+k-2},
\end{equation}
see \cite{Sta}.

If we think of the random walks as the motion of particles, then at 
each time we have a certain particle configuration. By considering these 
particles at all times we get a discrete, finite point process. The purpose 
of this paper is to show that this is a determinantal point process and 
compute the correlation kernel in terms of the associated Hahn polynomials,
\cite{JoNIP}, \cite{BKMcMi}. The derivation is based on the general framework
of \cite{JoDet} and a variant of the orthogonal polynomial method.
The main result is theorem \ref{thm3.1} below. 
The proof of that theorem also gives a proof of MacMahon's formula.
A certain continuous
scaling limit of this model, namely $a$ fixed and $b=c\to\infty$, converges
to a model of non-intersecting Brownian motions all started at the origin
and conditioned to end at the origin at time $T$. This Brownian motion
model is a transformation of Dyson's Hermitian Brownian motion model.
We will discuss these models in the next section and indicate how the
correlation kernel can be computed in these models using Hermite
polynomials and the orthogonal polynomial method. The result in this case
is closely related to the work in \cite{EyMe}, see also \cite{FNH}.
In the last section we will consider the discrete model where the
orthogonal polynomial method is less obvious. At the end of that section we 
will give some remarks concerning asymptotics.

\section{General framework and Dyson's Brownian motion}
\subsection{General framework}

Let $X_r$, $0\le r\le m$ be subsets of $\mathbb{R}$, $\phi_{r,r+1}:
X_r\to X_{r+1}$, $0\le r<m$, given functions and $\mu_r$ a measure on 
$X_r$, $1\le r\le m$, e.g. Lebesgue or counting measure. An element
$\underline{x}=(x^1,\dots,x^{m-1})\in X_1^n\times X_2^n\times\cdots
\times X_{m-1}^n\doteq \mathcal{X}$ is called a configuration. We think of
$x^r_1,\dots x^r_n$, $x^r=(x^r_1,\dots,x^r_n)$, as the positions of
particles in $X_r$, which we will call line $r$. 
Let $x^0\in X_0^n$ and $x^m\in X_m^n$ be fixed configurations, the initial and
final configurations respectively.
Define $\phi_{r,s}:X_r\times X_s\to\mathbb{R}$ for $r<s$ by
\begin{equation}\label{2.1}
\phi_{r,s}(x,y)=\int\phi_{r,r+1}(x,z_1)\dots\phi_{s-1,s}(z_{r-s-1},y)
d\mu_{r+1}(z_1)\dots d\mu_{s-1}(z_{r-s-1}),
\end{equation}
and $\phi_{r,s}\equiv 0$ if $r\ge s$. We will consider probability
measures on $\mathcal{X}$ of the form
\begin{equation}\label{2.2}
\frac 1{Z_{n,m}}\prod_{r=0}^{m-1}
\det(\phi_{r,r+1}(x^r_i,x^{r+1}_j))_{i,j=1}^nd\mu_1^n(x^1)\dots d\mu_{m-1}^n(
x^{m-1}),
\end{equation}
where $Z_{n,m}$ is a normalization constant. It is proved in \cite{JoDet}
that the measure (\ref{2.2}) has determinantal correlation functions,
i.e. the probability density with respect to the reference measure
$d\mu_{r_1}(y_1)\dots d\mu_{r_k}(y_k)$ of finding particles at
$z_1=(r_1,y_1),\dots, z_k=(r_k,y_k)$ is given by
\begin{equation}\label{2.3}
\det(K_{n,m}(z_i;z_j))_{i,j=1}^k
\end{equation}
where $K$ is the so called correlation kernel. This kernel is given by
\begin{equation}\label{2.4}
K_{n,m}(r,x;s,y)=-\phi_{r,s}(x,y)+\sum_{i,j=1}^n\phi_{r,m}(x,x_i^m)
(A^{-1})_{i,j}\phi_{0,s}(x_j^0,y),
\end{equation}
where $A=(\phi_{0,m}(x^0_i,x^m_j))_{i,j=1}^n$. Note that the kernel is 
not unique. We can multiply it by $\psi(r,x)/\psi(s,y)$ for an arbitrary
function $\psi\neq 0$ and get the same correlation functions.

\subsection{Dyson's Hermitian Brownian motion}
Let $H(t)$ be an $n\times n$ Hermitian matrix whose elements evolve
according to indepenent Ornstein-Uhlenbeck processes, see \cite{Dy},
\cite{Me}. We consider the stationary case. The probability measure
for seeing the matrices $H_1$,\dots, $H_{m-1}$ at times $t_1<\dots <t_m$
is 
\begin{equation}\label{2.5}
\frac 1{Z_{n,m}} e^{-\text{tr\,} H_1^2}\prod_{j=1}^{m-2}
\exp\left(-\frac{\text{tr\,}
(H_{j+1}-q_jH_j)^2}{1-q_j^2}\right)dH_1\dots dH_{m-1},
\end{equation}
where $dH_j$ is the Lebesgue measure on the space of Hermitian matrices,
and $q_j=\exp(-(t_{j+1}-t_j))$, $1\le j\le m-2$. Integrating out the angular
variables using the HarishChandra/Itzykson-Zuber formula, \cite{Me},
gives the eigenvalue measure
\begin{equation}\label{2.6}
\frac 1{Z_{n,m}'} \Delta_n(\lambda^1)\prod_{j=1}^n e^{-(\lambda^1_j)^2}
\prod_{r=1}^{m-2}\det \left(\exp\left(-\frac{(\lambda^{r+1}_j-q_r\lambda^r_i)
^2}{1-q_r^2}\right)\right)_{i,j=1}^n\Delta_n(\lambda^{m-1}),
\end{equation}
where $\Delta_n(\lambda)=\prod_{1\le i<j\le n}(\lambda_i-\lambda_j)$ is 
the Vandermonde determinant, and $\lambda_j^r$, $1\le j\le n$, are the 
eigenvalues of $H_r$.

If we set $\phi_{0,1}(i,x)=p_i(x)e^{-x^2}$, $\phi_{m-1,m}(x,i)=p_i(x)$,
where $p_i$ is a polynomial of degree $i$, $\phi_{r,r+1}(x,y)=
(\pi (1-q_r^2))^{-1/2}\exp(-(y-q_rx)^2/(1-q_r^2))$, $X_0=X_m=
\{0,\dots, n-1\}$, $x_1^0=x_i^m=i-1$, $X_r=\mathbb{R}$, $1\le r<m$ and
$\mu_r$ the Lebesgue measure, we see that (\ref{2.6}) is of the form
(\ref{2.2}). This is a basic example of a measure of the form
(\ref{2.2}). Here we have used the classical trick in the orthogonal
polynomial method in random matrix theory to write the Vandermonde
determinant as $\Delta_n(\lambda)=\det(p_i(\lambda_j))$. The polynomials
$p_i$ can be arbitrary but we choose them to be the normalized Hermite
polynomials. This will lead to a formula for the kernel (\ref{2.4}) in
terms of the Hermite polynomials. The key is the expansion, see e.g.
\cite{AAR},
\begin{equation}\label{2.7}
\frac 1{\sqrt{\pi(1-q^2)}}e^{-\frac{(qx-y)^2}{1-q^2}}=\sum_{k=0}^\infty
p_k(x)p_k(y)q^ke^{-y^2},
\end{equation}
$0<q<1$. Repeated use of this identity gives
\begin{equation}\label{2.8}
\phi_{0,s}(j,y)=e^{-j(t_s-t_1)}p_j(y)e^{-y^2}.
\end{equation}
Similarly,
\begin{equation}\label{2.9}
\phi_{r,m}(x,j)=e^{-j(t_{m-r}-t_r)}p_j(x).
\end{equation}
Using the orthonormality we obtain $\phi_{0,m}(i,j)=\exp(-j(t_{m-1}-t_1))
\delta_{ij}$ and hence $(A^{-1})_{ij}=\exp(j(t_{m-1}-t_1))\delta_{ij}$.
It also follows from (\ref{2.7}) that if $1\le r<s<m$, then
\begin{align}\label{2.10}
\phi_{r,s}(x,y)&=\frac 1{\sqrt{\pi(1-e^{2(t_r-t_s)})}}\exp\left(-\frac{(
e^{t_r-t_s}x-y)^2}{1-e^{2(t_r-t_s)}}\right)\notag\\
&=\sum_{k=0}^\infty p_k(x)p_k(y)e^{k(t_r-t_s)}e^{-y^2}.
\end{align}
Set $\chi_{t,s}=1$ if $t<s$ and  $\chi_{t,s}=0$ if $t\ge s$. From (\ref{2.4})
we get the {\it extended Hermite kernel},
\begin{align}\label{2.11}
K_{\text{ext.Herm.}}(t,x;s,y)&=-
\frac 1{\sqrt{\pi(1-e^{2(t-s)})}}\exp\left(-\frac{(
e^{t-s}x-y)^2}{1-e^{2(t-s)}}\right)\chi_{t,s}\notag\\
&+\sum_{k=0}^{n-1} e^{k(t-s)}p_k(x)p_k(y)e^{-y^2}
\end{align}
Using the second equality in (\ref{2.10}) we obtain the alternative
formula
\begin{equation}\label{2.11´}
K_{\text{ext.Herm.}}(t,x;s,y)=
\begin{cases} 
\sum_{k=0}^{n-1} e^{k(t-s)}p_k(x)p_k(y) e^{-y^2} &,t\ge s\\
-\sum_{k=n}^\infty  e^{k(t-s)}p_k(x)p_k(y) e^{-y^2} &,t<s.
\end{cases}
\end{equation}
Multiplying with $\exp(-x^2/2+y^2/2)$ we get the ordinary Hermite kernel
when $t=s$.

Let $\gamma$ be a positively oriented circle around the origin with 
radius $r>0$, and $\Gamma$ the line $\mathbb{R}\ni\to L+it$ with $L>r$.
Using the integral formulas, see \cite{AAR},
\begin{equation}
H_n(x)=\frac{2^n}{i\sqrt{\pi}} e^{x^2}\int_\Gamma e^{w^2-2xw}w^ndw,
\notag
\end{equation}
\begin{equation}
H_n(x)=\frac{n!}{2\pi i}\int_\gamma e^{-z^2+2xz}\frac{dz}{z^{n+1}},
\notag
\end{equation}
and $p_n(x)=(\sqrt{\pi}2^n n!)^{-1/2}H_n(x)$, it is not difficult to show that
\begin{align}\label{2.12}
&K_{\text{ext.Herm.}}(t,x;s,y)=-
\frac 1{\sqrt{\pi(1-e^{2(t-s)})}}\exp\left(-\frac{(
e^{t-s}x-y)^2}{1-e^{2(t-s)}}\right)\chi_{t,s}\notag\\
&+\frac{2}{(2\pi i)^2}\int_\Gamma dw\int_\gamma dz\frac{w^n}{z^n}
\frac 1{w-z}e^{w^2-2yw-e^{2(t-s)}z^2+2e^{t-s}xz}.
\end{align}
This double contour integral can be useful for asymptotic computations, 
for example to show convergence to the extended Airy kernel when
we have the edge scaling. To our knowledge the details for this has not 
been presented in the litterature, but using (\ref{2.12}) and the integral
formula for the extended Airy kernel it should be possible to do
this similarly to what was done for the extended Airy kernel in
\cite{JoAzt}.

\subsection{Non-intersecting Brownian motions}
A second closely related example is the following which involves 
non-intersecting Brownian motions. Consider $n$ non-intersecting Brownian
motions started at $x^0_i=\epsilon (i-1)$, $1\le i\le n$, at time 0 and 
conditioned to end at the same points at time $T$. Let $x^r_i$, 
$1\le i\le n$, 
denote the positions at time $\tau_r$, $1\le r<m$, where
$0=\tau_0<\tau_1<\dots<\tau_{m-1}<\tau_m=T$. By the Karlin-McGregor theorem
the probability density for $\underline{x}=(x^1,\dots,x^{m-1})\in(\mathbb{R}^n)
^{m-1}$ is given by
\begin{equation}\label{2.13}
\frac 1{Z^{\epsilon}_{n,m}}\prod_{r=0}^{m-1}\det(p_{\tau_{r+1}-\tau_r}
(x^r_i,x^{r+1}_j))_{i,j=1}^n.
\end{equation}
In the limit $\epsilon\to 0+$, corresponding to all particles starting
at the origin at time 0 and ending at the origin at time $T$, we get the 
probability density
\begin{align}\label{2.14}
\frac 1{Z_{n,m}}&\Delta_n(y^1)\prod_{j=0}^{n-1}e^{-(y^1_j)^2/2\tau_1}
\prod_{r=1}^{m-2}\det(e^{-(y^{r+1}_j-y^r_i)^2/2(\tau_{r+1}-\tau_r)})_{i,j=1}^n
\notag\\
&\times\Delta_n(y^{m-1})\prod_{j=0}^{n-1}e^{-(y^{m-1}_j)^2/2(T-\tau_{m-1})}.
\end{align}
This has again the general form (\ref{2.2}) with 
$\phi_{0,1}(i,y)=q_i(y)\exp(-y^2/2\tau_1)$,
$\phi_{r,r+1}(x,y)=\exp(-(y-x)^2/2(\tau_{r+1}-\tau_r))$ and 
$\phi_{m-1,m}(x,i)=\tilde{q}_i(x)\exp(-x^2/2(T-\tau_{m-1}))$, where
$q_i$ and $\tilde{q}_i$ are polynomials of degree $i$. 
The measure (\ref{2.14})
is actually a transformation of the measure (\ref{2.6}). Define
\begin{equation}\label{2.15}
d_r=\sqrt{\frac T{2\tau_r(T-\tau_r)}},
\end{equation}
$1\le r<m$, and $\tau_r=T(1+e^{-2t_r})^{-1}$. If we set $\lambda^r_j=
y^r_jd_r$, $1\le j\le n$, $1\le r<m$, then a straightforward computation
shows that (\ref{2.14}) transforms into (\ref{2.6}). In this way we
can also transform the extended Hermite kernel (\ref{2.11}) into a correlation
kernel for (\ref{2.14}). However, let us indicate how we can obtain it 
directly.

Set
\begin{equation}
c_{r,j}=\pi^{1/2}\left(\frac{\tau_r(T-\tau_{r+1})}{\tau_{r+1}(T-\tau_r)}
\right)^{j/2}.
\notag
\end{equation}
Then
\begin{equation}\label{2.16}
\int_{\mathbb{R}}e^{-\frac{(y-x)^2}{2(\tau_{r+1}-\tau_r)}}e^{-\frac{x^2}
{2\tau_r}}p_j(xd_r)dx=
c_{r,j}p_j(yd_{r+1})e^{-\frac{y^2}
{2\tau_{r+1}}},
\end{equation}
where $p_j$ is the $j$:th normalized Hermite polynomial.
This can be deduced from the identity
\begin{equation}\label{2.17}
\int_{\mathbb{R}}e^{-(x-y)^2}p_n(\alpha x)dx=
\pi^{1/2}(1-\alpha^2)^{n/2} p_n(\frac{\alpha y}{(1-\alpha^2)^{1/2}}),
\end{equation}
which in turn follows easily from the generating function for the Hermite
polynomials. Choose $q_j(x)=p_j(xd_1)$ and $\tilde{q}_j(x)=p_j(xd_{m-1})$.
It follows from (\ref{2.16}) that
\begin{equation}\label{2.18}
\phi_{0,s}(j,x)=\left(\frac{\tau_1}{\tau_s}\prod_{i=1}^{s-1}
(\tau_{i+1}-\tau_i)\right)^{1/2}2^{(s-1)/2}
\prod_{i=1}^{s-1} c_{i,j}p_j(xd_s)e^{-x^2/2\tau_s}
\end{equation}
and
\begin{equation}\label{2.19}
\phi_{r,m}(x,j)=\left(\frac{T-\tau_{m-1}}{T-\tau_r}\prod_{i=r}^{m-2}
(\tau_{i+1}-\tau_i)\right)^{1/2}2^{(m-r-1)/2}
\prod_{i=r}^{m-2} c_{i,j}p_j(xd_r)e^{-x^2/2(T-\tau_r)}.
\end{equation}
Using the orthogonality of the $p_j$:s and the general formula
(\ref{2.4}) we obtain the following expression for the correlation kernel
\begin{align}\label{2.20}
&K_{\text{BM}}(\tau_r,x;\tau_s,y)=-\frac 1{\sqrt{2\pi(\tau_r-\tau_s)}}
e^{-\frac{(x-y)^2}{2(\tau_r-\tau_s)}}\chi_{\tau,\tau_s}
\notag\\
&+\sum_{j=0}^{n-1}\left(\frac{\tau_r(T-\tau_s)}{\tau_s(T-\tau_r)}\right)^{j/2}
\left(\frac T{2\tau_s(T-\tau_r)}\right)^{1/2}
p_j(xd_r)p_j(yd_s)e^{-x^2/2(T-\tau_r)-y^2/2\tau_s}.
\end{align}
Here we have multiplied by the unimportant factor
$$
(2\pi)^{\frac{r-s}2}\left(\frac{(\tau_s-\tau_{s-1})\dots (\tau_2-\tau_1)}
{(\tau_r-\tau_{r-1})\dots (\tau_2-\tau_1)}\right)^{1/2}.
$$

If we go back to the transformation discussed above we see that
\begin{equation}\label{2.21}
\frac 1{\sqrt{d_rd_s}}K_{\text{BM}}(\tau_r,\frac x{d_r};\tau_s,\frac y{d_s})
e^{\frac{x^2\tau_r}T-\frac{y^2\tau_s}T}\left(\frac{T/\tau_r-1}
{T/\tau_s-1}\right)^{1/4}=K_{\text{ext.Herm.}}(t_r,x;t_s,y),
\end{equation}
with $\tau_r=T(1+\exp(-2t_r))^{-1}$ and $K_{\text{ext.Herm.}}$ given
by (\ref{2.11}).

\section{The extended Hahn kernel}
\subsection{Derivation of the kernel}
Consider $a$ symmetric, simple random walks with initial points $(0,2j)$ and 
final points $(b+c, c-b+2j)$, $0\le j\le a-1$, conditoned not to intersect
in the whole time interval $[0,b+c]$. The single step transition kernel for 
one particle is
\begin{equation}\label{3.1}
\frac 12\phi(x,y)=\frac 12\delta_{x-1,y}+\frac 12\delta_{x+1,y}.
\end{equation}
The configuration at time $t=r$, which we also call the configuration
on the $r$:th line, is given by points $z^r_j$, $0\le j<a$,
$z^r_0<\dots< z^r_{a-1}$, where $z^0_j=2j$, $z^{b+c}_j=c-b+2j$.
We think of these points as the positions of particles.
By the Lindstr\"om-Gessel-Viennot method, \cite{Ste}, our probability
measure on the set of configurations $\underline{z}=(z^r_j)$ in
$(\mathbb{Z}^a)^{b+c-1}$ is
\begin{equation}\label{3.2}
p(\underline{z})=\frac 1{Z(a,b,c)}\prod_{r=0}^{b+c-1}\det(\phi
(z^r_j,z^{r+1}_k))_{j,k=0}^{a-1}.
\end{equation}
Here $Z(a,b,c)$ is the total number of configurations and is given by 
MacMahon's formula (\ref{1.1}).

The measure (\ref{3.2}) has exactly the general form (\ref{2.2}) (with $\mu$
counting measure on $\mathbb{Z}$), and we want to
compute the correlation kernel (\ref{2.4}). To do this we will use the 
orthogonal polynomial method in a similar way that was used for the
non-intersecting Brownian motions in the last section. How this should be done
is not obvious from (\ref{3.2}). It is shown in \cite{JoDOPE}, that the
induced probability ensemble on a single line is an orthogonal polynomial 
ensemble, where the relevant polynomials are the associated Hahn polynomials.
This indicates that we should modify the first and the last
factors in (\ref{3.2}) by doing row operations so that we get a situation 
where
the matrix $A$ in (\ref{2.4}) is diagonal.

The normalized associated Hahn polynomials, \cite{Niki}, \cite{JoNIP},
\cite{BKMcMi}, can be defined using a hypergeometric function by
\begin{align}\label{3.3}
q^{(\alpha,\beta)}_{n,N}(x)&=
\frac{(-N-\beta)_n(-N)_n}{d^{(\alpha,\beta)}_{n,N}
n!}{}_3F_2 
\left(\begin{matrix} -n,n-2N-\alpha-\beta-1,-x \\
-N-\beta,-N \end{matrix}\,;\, 1\right)
\notag\\
&=\frac{(-N-\beta)_n(-N)_n}{d^{(\alpha,\beta)}_{n,N}n!}
\sum_{j=0}^n\binom{n}{j}(-1)^j\frac{(-x)_j(n-2N-\alpha-\beta-1)_j}
{(-N-\beta)_j(-N)_j},
\end{align}
where
\begin{equation}\label{3.4}
\left(d^{(\alpha,\beta)}_{n,N}\right)^2=\frac{(\alpha+\beta+N+1-n)_{N+1}}
{(\alpha+\beta+2N+1-2n)n!(\beta+N-n)!(\alpha+N-n)!(N-n)!},
\end{equation}
and we use the standard notation $(a)_n=a(a+1)\dots(a+n-1)$. These 
polynomials are orthogonal with respect to the weight
\begin{equation}\label{3.5}
w^{(\alpha,\beta)}_N(x)=\frac 1{x!(x+\alpha)!(N+\beta-x)!(N-x)!},
\end{equation}
on $\{0,1,\dots, N\}$, i.e.
\begin{equation}\label{3.6}
\sum_{x=0}^Nq^{(\alpha,\beta)}_{n,N}(x)q^{(\alpha,\beta)}_{m,N}(x)
w^{(\alpha,\beta)}_N(x)=\delta_{n,m},
\end{equation}
for $0\le n,m\le N$. Below we will sometimes use the convention that $1/n!=0$
if $n<0$, so that the summation in (\ref{3.6}) for example could be
extended to $x\in\mathbb{Z}$.

Our goal is to give a formula for the correlation kernel in terms of the 
associated Hahn polynomials. First, we need some notation. Let $a,b,c
\in\mathbb{Z}^+$, $b\le c$. Set $a_r=|c-r|$, $b_r=|b-r|$,
\begin{equation}\label{3.7}
\alpha_r=\begin{cases} -r &,0\le r\le b \\
r-2b &,b\le r\le b+c \end{cases}
\end{equation}
and
\begin{equation}\label{3.8}
\gamma_r=\begin{cases} r+a-1 &,0\le r\le b \\
b+a-1 &,b\le r\le c \\a+b+c-1-r &,c\le r\le b+c.
\end{cases}
\end{equation}
Define
\begin{equation}\label{3.9}
\omega_r(x)=\begin{cases} ((b_r+x)!(\gamma_r+a_r-x)!)^{-1} &,0\le r\le b \\
(x!(\gamma_r+a_r-x)!)^{-1}&,b\le r\le c \\
(x!(\gamma_r-x)!)^{-1} &,c\le r\le b+c
\end{cases}
\end{equation}
and
\begin{equation}\label{3.10}
\tilde{\omega}_s(x)=\begin{cases} 
(y!(\gamma_s-y)!)^{-1} &,0\le s\le b \\
((b_s+y)!(\gamma_s-y)!)^{-1} &,b\le s\le c \\
((b_s+y)!(\gamma_s+a_s-y)!)^{-1}&,c\le s\le b+c.
\end{cases}
\end{equation}

\begin{theorem}\label{thm3.1}
The point process on $(\mathbb{Z}^a)^{b+c-1}$ defined by (\ref{3.20}) has
determinantal correlation functions with kernel given by
\begin{align}\label{3.11}
&K_H(r,\alpha_r+2x;s,\alpha_s+2y)=-\phi_{r,s}(\alpha_r+2x,\alpha_s+2y)
\notag\\
&+\sum_{n=0}^{a-1}\sqrt{\frac{(a+s-1-n)!(a+b+c-r-1-n)!}{(a+r-1-n)!(a+b+c-1-n)!}
} q_{n,\gamma_r}^{(b_r,a_r)}(x)q_{n,\gamma_s}^{(b_s,a_s)}(y)\omega_r(x)
\tilde{\omega}_s(y),
\end{align}
for $0<r,s<b+c$, $x,y\in\mathbb{Z}$. Here $\phi_{r,s}\equiv 0$ if
$r\ge s$ and
\begin{equation}\label{3.11´}
\phi_{r,s}(x,y)=\binom{s-r}{\frac{y-x+s-r}2}
\end{equation}
if $r<s$.
\end{theorem}

\begin{proof}
Set
\begin{equation}
c_{j,k}=\frac 1{(a-k)(j-k)!(a-1-j)!},
\notag
\end{equation}
for $0\le j,j<a$,
\begin{equation}
f_{n,k}=\binom{n}{k}\frac{(n-2a-b-c+1)_k}{(-a-c+1)_k(-a)_k}
\notag
\end{equation}
and
\begin{equation}
f^\ast_{n,k}=\binom{n}{k}\frac{(n-2a-b-c+1)_k}{(-a-b+1)_k(-a)_k}.
\notag
\end{equation}
for $0\le k\le n$.
Define
\begin{align}\label{3.12}
\psi(n,z)&=\sum_{m=0}^n f_{n,m}\sum_{j=m}^{a-1} c_{j,m}\phi(2j,z),
\notag\\
\psi^\ast(n,z)&=\sum_{m=0}^n f^\ast _{n,m}\sum_{j=m}^{a-1} c_{j,m},
\phi(c-b++2j,z)
\end{align}
$0\le n<a$, $z\in\mathbb{Z}$.

We will now do row operations to modify the first and the last factor in
(\ref{3.2}).
\begin{align}\label{3.13}
\det(\phi(x^0_j,x^1_k))_{j,k=0}^{a-1}&=
\det(\phi(2m,x^1_k))_{m,k=0}^{a-1}
\notag\\
&=\det(\frac 1{c_{m,m}}\sum_{j=m}^{a-1}c_{j,m}\phi(2j,x^1_k))_{m,k=0}^{a-1}
\notag\\
&=\prod_{m=0}^{a-1}\frac 1{c_{m,m}}
\det(\sum_{j=n}^{a-1}c_{j,n}\phi(2j,x^1_k))_{n,k=0}^{a-1}
\notag\\
&=\prod_{m=0}^{a-1}\frac 1{c_{m,m}}
\det(\frac 1{f_{n,n}}\sum_{m=0}^{n}f_{n,m}
\sum_{j=m}^{a-1}c_{j,m}\phi(2j,x^1_k))_{n,k=0}^{a-1}
\notag\\
&=\prod_{m=0}^{a-1}\frac 1{c_{m,m}f_{m,m}}
\det(\psi(n,x^1_k))_{n,k=0}^{a-1}.
\end{align}
In the same way we obtain
\begin{equation}\label{3.14}
\det(\phi(x^{b+c-1}_j,x^{b+c}_k))_{j,k=0}^{a-1}
=\prod_{m=0}^{a-1}\frac 1{c_{m,m}f^\ast_{m,m}}
\det(\psi(n,x^{b+c-1}_k))_{n,k=0}^{a-1}.
\end{equation}
If we now set $\phi_{0,1}(n,y)=\psi(n,y)$, $\phi_{b+c-1,b+c}(y,n)=\psi^\ast
(n,y)$ and $\phi_{r,r+1}(x,y)=\phi(x,y)$, $1\le r<b+c-1$, the probability 
measure (\ref{3.2}) can be written
\begin{equation}\label{3.15}
p(\underline{y})=\frac 1{Z(a,b,c)}\prod_{m=0}^{a-1}\frac 1{c_{m,m}^2
f_{m,m}f^\ast_{m,m}}\prod_{r=0}^{b+c-1}\det(\phi_{r,r+1}(y_j^r,y_k^{r+1}))
_{j,k=0}^{a-1},
\end{equation}
where $y^0_j=y^{b+c}_j=j$, $0\le j<a$.

Write $\phi^{\ast n}(x,y)=\phi\ast\dots\ast\phi(x,y)$ ($n$ factors) if
$n\ge 2$, $\phi^{\ast 1}(x,y)=\phi(x,y)$ and $\phi^{\ast 0}(x,y)=
\delta_{x,y}$.
We want to compute $\phi_{0,s}$, $\phi_{r,b+c}$ and $\phi_{0,b+c}$
for $1\le r,s<b+c$. By definition
\begin{equation}\label{3.16}
\phi_{0,r}(n,y)=\sum_{z\in\mathbb{Z}}\psi(n,z)\phi^{\ast(r-1)}(z,y),
\end{equation}
and
\begin{equation}\label{3.17}
\phi_{r,b+c}(y,n)=\sum_{z\in\mathbb{Z}}\psi^\ast(n,z)\phi^{\ast(b+c-r-1)}(z,y),
\end{equation}
since $\phi(x,y)=\phi(y,x)$.

\begin{claim}\label{cl3.2}
If $z\in 2\mathbb{Z}+1$, then
\begin{equation}\label{3.18}
\psi(n,z)=\sum_{j=0}^n \binom{n}{j}\frac{(n-2a-b-c+1)_j}
{(-a-c+1)_j(-a)_j(\frac{z+1}2-j)!(a-\frac {z+1}2)!},
\end{equation}
and if $z\in 2\mathbb{Z}$, then $\psi(n,z)=0$.
\end{claim}

\begin{proof}
By definition
\begin{equation}
\psi(n,z)=\sum_{m=0}^n\binom{n}{m}\frac{(n-2a-b-c+1)_j}{(-a-c-1)_m(-a)_m}
\sum_{j=m}^{a-1}\frac{\phi(2j,z)}{(a-m)(j-m)!(a-1-j)!}.
\notag
\end{equation}
Now, with $z=2\zeta-1$,
\begin{align}
&\sum_{j=m}^{a-1}\frac{\delta_{2j-1,z}+\delta_{2j+1,z}}{(a-m)(j-m)!(a-1-j)!}
\notag\\
&=\frac 1{(a-m)(\zeta-m)!(a-1-\zeta)!}+
\frac 1{(a-m)(\zeta-1-m)!(a-\zeta)!}
\notag\\
&=\frac 1{(\frac{z+1}2-m)!
(a-\frac{z+1}2)!},
\notag
\end{align}
and (\ref{3.18}) follows.
\end{proof}

If $z\in 2\mathbb{Z}-1$, then simailarly
\begin{equation}\label{3.19}
\psi^\ast(n,c-b+z)=\sum_{j=0}^n \binom{n}{j}\frac{(n-2a-b-c+1)_j}
{(-a-b+1)_j(-a)_j(\frac{z+1}2-j)!(a-\frac {z+1}2)!}.
\end{equation}

\begin{claim}\label{cl2}
\begin{equation}\label{3.20}
\phi_{0,r}(n,y)=(a+1)_{r-1}\sum_{j=0}^n\binom{n}{j}
\frac{(n-2a-b-c+1)_j}{(-a-c+1)_j(-a-r+1)_j}
\frac 1{(\frac{y+r}2-j)!(a-1-\frac{y-r}2)!}.
\end{equation}
\end{claim}

\begin{proof}
By induction on $r$. The statement is true for $r=1$ by (\ref{3.18}). We have
\begin{align}
&\phi_{0,r+1}(n,y)=\sum_{x\in\mathbb{Z}}\phi_{0,r}(n,x)\phi(x,y)
\notag\\
&\sum_{x\in\mathbb{Z}}\phi_{0,r}(n,x)(\delta_{x,y+1}+\delta_{x,y-1})
=\phi_{0,r}(n,y+1)+\phi_{0,r}(n,y-1)
\notag\\
&=(a+1)_{r-1}\sum_{j=0}^n\binom{n}{j}\frac{(n-2a-b-c+1)_j}
{(-a-c+1)_j(-a-r+1)_j}\frac 1{(\frac{y+r+1}2-j)!(a-1-\frac{y-r-1}2)!}
\notag\\
&\times [a-\frac{y-r+1}2+\frac{y+r+1}2-j].
\notag
\end{align}
Now,
\begin{equation}
\frac{(a+1)_{r-1}}{(-a-r+1)_j}(a-r-j)=\frac{(a+1)_r}{(-a-r)_j},
\notag
\end{equation}
and the claim is proved.
\end{proof}

Also,
\begin{align}
&\phi_{r,b+c}(c-b+x,n)=\sum_{z\in\mathbb{Z}}
\psi^\ast(n,z)\phi^{\ast(b+c-r-1)}(z,c-b+x)
\notag\\
&=\sum_{z\in\mathbb{Z}}\psi^\ast(n,c-b+z))\phi^{\ast(b+c-r-1)}
(c-b+z,c-b+x)
\notag\\
&==\sum_{z\in\mathbb{Z}}\psi^\ast(n,c-b+z))\phi^{\ast(b+c-r-1)}(z,x).
\notag
\end{align}
We can now proceed exactly as in the proof of claim \ref{cl2} and show
that
\begin{align}\label{3.21}
\phi_{r,b+c}(y,n)&=(a+1)_{b+c-r-1}\sum_{j=0}^n\binom{n}{j}
\frac{(n-2a-b-c+1)_j}{(-a-b+1)_j(-a-b-c+r+1)_j}
\notag\\
&\times\frac 1{(\frac{y-r}2+b-j)!(a+c-1-\frac{y+r}2)!}.
\end{align}

Introduce new coordinates, which we will call the {\it Hahn coordinates}
on line $r$ by
\begin{equation}
x^r_k=\frac{y^r_k-\alpha_r}2.
\notag
\end{equation}
Then, $0\le x^r_k\le\gamma_r$. One motvation to use these coordinates is that
it is easier to recognize the Hahn polynomials when using them.
Since $\phi_{0,r}(i,z)$ is zero unless $z+r$ is even, i.e. unless
$z-\alpha_r$ is even, we obtain
\begin{equation}\label{3.22}
A_{nm}=\sum_{z\in\mathbb{Z}}\phi_{0,r}(n,\alpha_r+2z)\phi_{r,b+c}
(\alpha_r+2z,m).
\end{equation}
The correlation kernel is given by
\begin{align}\label{3.23}
&K(r,2x+\alpha_r;s,2y+\alpha_s)=-\phi_{r,s}(2x+\alpha_r,2y+\alpha_s)
\notag\\
&+\sum_{i,j=0}^{a-1}\phi_{r,b+c}(2x+\alpha_r,i)(A^{-1})_{ij}
\phi_{0,s}(j,2y+\alpha_s)
\end{align}
according to (\ref{2.4}). We want to express $\phi_{0,r}(j,2y+\alpha_r)$
and $\phi_{r,b+c}(2x+\alpha_r,i)$ in terms of the associated 
Hahn polynomials. 
In order to do so we have to distinguish three cases, $1\le r\le b$, $b\le r
\le c$ and $c\le r\le b+c$.

Set $a_r=|c-r|$ and $b_r=|b-r|$.

\noindent
(i) Consider first the case $1\le r\le b$. By (\ref{3.3}) and
(\ref{3.20})
\begin{align}\label{3.24}
&\phi_{0,r}(n,\alpha_r+2z)=(a+1)_{r-1}\sum_{j=0}^n\binom{n}{j}
\frac{(n-2a-b-c+1)_j}{(-a-c+1)_j(-a-r+1)_j}
\notag\\
&\times\frac 1{(z-j)!(a+r-1-z)!}=
\frac{(a+1)_{r-1}d_{n,\gamma_r}^{(b_r,a_r)}n!}{(-a-c+1)_n(-a-r+1)_n}
q_{n,\gamma_r}^{(b_r,a_r)}(z)\frac 1{z!(\gamma_r-z)!}.
\end{align}
Also, by (\ref{3.21}),
\begin{align}
&\phi_{r,b+c}(\alpha_r+2z,n)=
(a+1)_{b+c-r-1}\sum_{j=0}^n\binom{n}{j}
\frac{(n-2a-b-c+1)_j}{(-a-b+1)_j(-a-b-c+r+1)_j}
\notag\\
&\times\frac 1{(b-r+z-j)!(a+c-1-z)!}
\notag\\
&=\frac{(a+1)_{b+c-r-1}}{(b-r+z)!(a+c-1-z)!}\sum_{j=0}^n
\frac{(-n)_j(n-2a-b-c)_j(-b+r-z)_j}
{j!(-a-b+1)_j(-a-b-c+r+1)_j}
\notag\\
&=\frac{(a+1)_{b+c-r-1}}{(b_r+z)!(\gamma_r+\alpha_r-z)!}
 {}_3F_2 \left(\begin{matrix} -n.n-2a-b-c+1,-b+r-z \\
-a-b+1,-a-b-c+r+1 \end{matrix}\,;\,1\right).
\notag
\end{align}
We can rewrite this using the following hypergeometric identity,
\cite{AAR} p. 141,
\begin{equation}\label{3.24´}
 {}_3F_2\left(\begin{matrix} -n,a,b \\d,e\end{matrix};1\right)
=\frac{(d-a)_n(e-a)_n}{(d)_n(e)_n}
 {}_3F_2\left(\begin{matrix} -n,a,a+b-n-d-e+1 
\\a-n-d+1,a-n-e+1\end{matrix};1\right).
\end{equation}
This gives
\begin{align}\label{3.25}
&\phi_{r,b+c}(\alpha_r+2z)=\frac{(a+1)_{b+c-r-1}(a+c-n)_n(a+r-n)_n}
{(-a+b+1)_n(-a-b-c+r+1)_n(b_r+z)!(\gamma_r+a_r-z)!}
\notag\\
&\times
 {}_3F_2\left(\begin{matrix} -n,n-2\gamma_r-a_r-b_r-1,-z 
\\-\gamma_r-a_r,-\gamma_r\end{matrix};1\right).
\notag\\
&=\frac{(a+1)_{b+c-r-1}(a+c-n)_n(a+r-n)_nd_{n,\gamma_r}^{(b_r,a_r)}n!}
{(-a+b+1)_n(-a-b-c+r+1)_n(-a-c+1)_n(-a-r+1)_n}
\notag\\
&\times q_{n,\gamma_r}^{(b_r,a_r)}(z)\frac 1{(b_r+z)!(\gamma_r+a_r-z)!}
\end{align}
We can now compute $A_{nm}$ given by (\ref{3.22}) by picking $r$ between
1 and $b$, the choice does not matter.
Using (\ref{3.6}), (\ref{3.24}) and (\ref{3.25}) we obtain, after
some simplification
\begin{equation}\label{3.26}
A_{nm}=C_n(a,b,c)^{-1}\delta_{n,m},
\end{equation}
where
\begin{equation}\label{3.27}
C_n(a,b,c)=\frac{(a+b-1)!(a+c-1)!(2a+b+c-2n-1)a!^2}{n!(2a+b+c-n-1)!}.
\end{equation}

\noindent
(ii) Next we consider the case $b\le r\le c$.
The computations are similar to those in the previous case. We find
\begin{align}\label{3.28}
\phi_{0,r}(n,\alpha_r+2z)&=
\frac{(a+1)_{r-1}(a+b-n)_n(a+b+c-r-n)_n d_{n,\gamma_r}^{(b_r,a_r)}n!}
{(-a-c+1)_n(-a-r+1)_n(-a-b+1)_n(-a-b-c+1+r)_n}
\notag\\
&\times q_{n,\gamma_r}^{(b_r,a_r)}(z)\frac 1{(b_r+z)!(\gamma_r-z)!}.
\end{align}
Here we have used the hypergeometric identity
(\ref{3.24´}). Also, we find
\begin{align}\label{3.29}
\phi_{r,b+c}(\alpha_r+2z,n)&=
\frac{(a+1)_{b+c-r-1} d_{n,\gamma_r}^{(b_r,a_r)}n!}
{(-a-b+1)_n(-a-b-c+1+r)_n}
\notag\\
&\times q_{n,\gamma_r}^{(b_r,a_r)}(z)\frac 1{z!(\gamma_r+a_r-z)!}.
\end{align}

\noindent
(iii) Finally we come to the case $c\le r\le b+c$, and again the computations
are similar. We obtain
\begin{align}\label{3.30}
\phi_{0,r}(n,\alpha_r+2z)&=
\frac{(a+1)_{r-1}(a+b-n)_n(a+b+c-r-n)_n d_{n,\gamma_r}^{(b_r,a_r)}n!}
{(-a-c+1)_n(-a-r+1)_n(-a-b+1)_n(-a-b-c+1+r)_n}
\notag\\
&\times q_{n,\gamma_r}^{(b_r,a_r)}(z)\frac 1{(b_r+z)!(\gamma_r+a_r-z)!},
\end{align}
where we have used the identity (\ref{3.24´}). Also,
\begin{align}\label{3.31}
\phi_{r,b+c}(\alpha_r+2z,n)&=
\frac{(a+1)_{b+c-r-1} d_{n,\gamma_r}^{(b_r,a_r)}n!}
{(-a-b+1)_n(-a-b-c+1+r)_n}
\notag\\
&\times q_{n,\gamma_r}^{(b_r,a_r)}(z)\frac 1{z!(\gamma_r-z)!}.
\end{align}

We now have all the ingredients in (\ref{2.4}). It follows from
(\ref{3.26}) that
\begin{equation}
(A^{-1})_{ij}=C_i(a,b,c)\delta_{ij},
\notag
\end{equation}
and some computation now gives (\ref{3.11}). Note that $\phi_{r,s}(x,y)$ is 
the number of random walk paths from $x$ to $y$ in $s-r$ steps and hence 
is given by (\ref{3.11´}).

The computations in the proof of the theorem also gives a proof of 
MacMahons formula. We have
\begin{equation}\label{3.31'}
Z(a,b,c)=\prod_{n=0}^{a-1}\frac 1{c_{n,n}^2d_{n,n}d^\ast_{n,n}}
\det A.
\end{equation}
A computation gives
\begin{equation}
\prod_{n=0}^{a-1}\frac 1{c_{n,n}^2d_{n,n}d^\ast_{n,n}}
=\prod_{n=0}^{a-1}\frac{(2a+b+c-2n-1)!^2(a+b-1)!(a+c-1)!}
{(a+b-1-n)!(a+c-1-n)!(2a+b+c-n-1)!^2}.
\notag
\end{equation}
It follows from (\ref{3.26}) and (\ref{3.27}) that
\begin{equation}
\det A
=\prod_{n=0}^{a-1}\frac{n!(2a+b+c-n-1)!}{(a+b-1)!(a+c-1)!(2a+b+c-2n-1)a!^2}.
\notag
\end{equation}
Hence, by (\ref{3.31'}) and after some simplification
\begin{equation}
Z(a,b,c)=\prod_{n=0}^{a-1}\frac{n!(b+c+n)!}{(b+n)!(c+n)!},
\notag
\end{equation}
which is the same as (\ref{1.1}).
\end{proof}.

\subsection{Some remarks about asymptotics}

As discussed above the non-intersecting Brownian motion model (\ref{2.14}) 
is a kind of continuum version  of the random walk model. In fact it can be
obtained as a scaling limit of the random walk model. For the associated 
Hahn polynomials we have the asymptotics
\begin{equation}\label{3.32}
\lim_{N\to\infty}d_{n,N}^{(\alpha,\alpha)} n!\left(
-\frac{2}{N^{3/2}\sqrt{(2t+1)(t+1)}}\right)^np_{n,N}^{(\alpha,\alpha)}
(\frac N2+2z\sqrt{\frac{2t+1}{t+1}N})=H_n(z),
\end{equation}
where $\alpha/N\to t\ge 0$, uniformly for $z$ in a compact subset of 
$\mathbb{C}$. Here $H_n(z)$ is the ordinary Hermite polynomial of degree
$n$. This can be proved by a slight modification of the argument in
\cite{Holm} based on the recurrence relation. Using (\ref{3.32}) and
standard asymptotics for the binomial coefficient it follows that
\begin{equation}\label{3.33}
2^{r-s}\sqrt{\frac{k}{2T}} K_H(r,x;s,y)\to K_{\text{BM}}(\tau,\xi;\sigma,y)
\end{equation}
as $k\to\infty$ if $r/k\to2\tau/T$, $s/k\to 2\sigma/T$,
$x/\sqrt{k}\to\xi\sqrt{2/T}$, $y/\sqrt{k}\to\eta\sqrt{2/T}$,
where $K_{\text{BM}}$ is given by (\ref{2.20}). So in this sense we have 
convergence  to the Brownian motion model. It should also be possible to 
prove this directly, i.e. that the measure (\ref{3.2}) converges, when 
rescaled as above, to the measure (\ref{2.14}), compare the arguments in
\cite{KatTan}.

A more interesting, and also much more difficult limit is to consider
the case when, $a$, $b$ and $c$ go to infinity with the same rate, 
say $a=b=c\to\infty$. In particular it is interesting to consider the 
fluctuations of the top (and bottom) curves which bound the so called frozen 
regions, \cite{CLP}, in the tiling. If we restrict to a single line, this
has been done recently by \cite{BKMcMi} using very precise asymptotics
for Hahn polynomials derived using Riemann-Hilbert techniques. This shows
for example that if $n=a=b=c$ then the last (first) particle fluctuates
like $n^{1/3}$ in the appropriate region
and that the fluctuations are given by the 
Tracy-Widom distribution, see \cite{BKMcMi} for details. If these asymptotic
results could be extended to the extended (associated) Hahn kernel,
(\ref{3.11}), it should be possible to prove the convergence of the boundary
curve of the frozen region to the Airy process, \cite{PrSp}, \cite{JoDet},
as has been done for some other tiling problems in \cite{FeSp} and 
\cite{JoAzt}.


\begin{thebibliography}{99}

\itemsep=\smallskipamount

\bibitem{AAR} G.E. Andrews, R. Askey, R. Roy, {\em Special Functions,
Encyclopedia of Mathematics and its applications 71,} Cambridge
University Press, Cambridge, 1999


\bibitem{BKMcMi} J. Baik, T. Kriecherbauer, K.D.T.-R MacLaughlin, P. Miller,
{\em Uniform asymptotics for polynomials orthogonal with respect to a 
general class of discrete weights and universality results for 
associated ensembles,} math.CA/0310278

\bibitem{CLP} H. Cohn, M. Larsen, J. Propp, {\em The shape of a typical
boxed plane partition,} New York J. of Math., {\bf 4}, (1998), 137 -
165

\bibitem{Dy} F. J. Dyson, {\em A Brownian-Motion Model for the eigenvalues
of a Random Matrix,} J. Math. Phys., {\bf 3} (1962), 1191 - 1198

\bibitem{EyMe} B. Eynard, M.L. Mehta, {\em Matrices coupled in a chain
    I: Eigenvalue correlations}, J. of Phys. A, {\bf 31} (1998), 4449
    - 4456

\bibitem{FeSp} P. L. Ferrari, H. Spohn, {\em Step fluctuations for a
faceted crystal}, J. Stat. Phys., {\bf 113} (2003), 1 - 46

\bibitem{FNH} P.J. Forrester, T. Nagao, G. Honner, {\em Correlations
    for the orthogonal-unitary and symplectic-unitary transitions at
    the soft and hard edges}, Nucl. Phys. B, {\bf 553} (1999), 601 - 643

\bibitem{Holm} K. Holm{\aa}ker, {\em On a discrete Rodrigues' formula and a 
second class of orthogonal Hahn polynomials,} Preprint, Department of 
Mathematics, Chalmers University of Technology, No. 1977-12

\bibitem{JoDOPE} K. Johansson, {\em Discrete orthogonal polynomial ensembles
and the Plancherel measure,} 
Annals of Math., {\bf 153} (2001), 259 - 296

\bibitem{JoNIP} K. Johansson, {\em Non-intersecting paths, random
    tilings and random matrices}, 
    Probab.Theory Relat. Fields, {\bf 123} (2002), 225--280

\bibitem{JoDet} K. Johansson, {\em Discrete polynuclear growth and
determinantal processes,} Commun. Math. Phys., {\bf 242} (2003), 277 - 329

\bibitem{JoAzt} K. Johansson, {The Arctic circle and the Airy process,}
math.PR/0306216, to appear in Ann. Probab.

\bibitem{KatTan} M. Katori, H. Tanemura, {\em Scaling limit of vicious 
walks and two-matrix model,} Phys. Rev. E (2002)

\bibitem{Ke} R. Kenyon, {\em Local statistics of lattice dimers,}
  Ann. Inst. H. Poincar\'e, Probabilit\'es et Statistiques, {\bf 33}
  (1997), 591 - 618

\bibitem{Me} M. L. Mehta, {\em Random Matrices,} 2nd ed., Academic Press, 
San Diego 1991

\bibitem{Niki} A. F. Nikiforov, S. K. Suslov, V. B. Uvarov, {\em Classical
Orthogonal Polynomials of a Discrete Variable,} Springer Series in
Computational Physics, Springer-Verlag, Berlin Heidelberg, 1991

\bibitem{PrSp} M. Pr\"ahofer, H. Spohn, {\em Scale invariance of the
    PNG droplet and the Airy process}, J. Stat. Phys., {\bf 108}
(2002), 1076--1106

\bibitem{Sta} R. P. Stanley, {\em Enumerative Combinatorics,} Vol. 2,
Cambridge University Press, 1999

\bibitem{Ste} J. R. Stembridge, {\em Nonintersecting Paths,
    Pfaffians, and Plane Partitions,} Adv. in Math., {\bf 83} (1990), 96 - 131




\end{thebibliography}
\end{document}